\documentclass[11pt,a4paper]{article}
\usepackage{amsfonts}
\usepackage{}

\usepackage{graphicx}
\usepackage{amsfonts,amsmath, amssymb}

\newtheorem{theo}{\textbf{\ \ \quad Theorem}}[section]

\newtheorem{lem}{\textbf{\ \ \quad Lemma}}[section]
\newtheorem{remark}{\textbf{\ \ \quad Remark}}[section]
\newtheorem{col}{\textbf{\ \ \quad Corollary}}[section]
\newtheorem{prop}{\textbf{\ \ \quad Proposition}}[section]
\newtheorem{defi}{\textbf{\ \ \quad Definition}}[section]

\newcommand{\lbl}[1]{\label{#1}}

\newcommand{\be}{\begin{equation}}
\newcommand{\ee}{\end{equation}}
\newcommand\bes{\begin{eqnarray}}
\newcommand\ees{\end{eqnarray}}
\newcommand{\bess}{\begin{eqnarray*}}
\newcommand{\eess}{\end{eqnarray*}}

\newcommand{\nm}{\nonumber}

\newcommand{\ds}{\displaystyle}

{\setlength\arraycolsep{2pt}}

\title{Impacts of Noise on a Class of Partial Differential Equations }

\author{Guangying Lv$^a$ and Jinqiao Duan$^b$\\
\\
\ \\
   {\small \it $^a$ Institute of Contemporary Mathematics, Henan University}\\
  {\small \it Kaifeng, Henan 475001, China}\\
  {\small \tt gylvmaths@henu.edu.cn}\\
  {\small \it $^b$ Department of Applied Mathematics, Illinois Institute of Technology}\\
  {\small \it Chicago, IL 60616 }\\
   {\small \tt duan@iit.edu }
}

\begin{document}
\maketitle

\medskip

\begin{abstract}
This paper is concerned with effects of noise on the solutions of
partial differential equations.
We first provide a sufficient condition to ensure the existence of a unique
positive solution for a class of stochastic partial differential equations. Then, we prove that noise could induce
singularities (finite time blow up of solutions). Finally, we show that
a stochastic Allen-Cahn equation does not have finite time singularities and the unique solution exists globally.



\end{abstract}

\baselineskip=15pt

\section{Introduction}
\setcounter{equation}{0}

Stochastic partial differential equations (SPDEs) are playing an increasingly
important role in modeling complex phenomena in physics, geophysics and biology. In recent years,
existence, uniqueness, stability, blow-up phenomenon, invariant measures and other
properties of the solutions to SPDEs have been extensively investigated  \cite{Cb,LR,L}. It is known that the existence and uniqueness of
global solutions to SPDEs can  be   established under appropriate conditions (\cite{Cb,DW}).

It is also  known   that    certain deterministic parabolic or hyperbolic  partial
differential equations (even with polynomial nonlinearity)  tend to develop
singularities in finite time \cite{F,SGKM}. These equations
only have solutions which are defined locally in time. For example, consider the following   equation
   \bes\left\{\begin{array}{lll}
   u_t-\Delta u=u^{1+\alpha}, \ \ \qquad t>0,&x\in  D,\\[1.5mm]
   u(x,0)=u_0(x), \ \ \ &x\in D,\\[1.5mm]
   u(x,t)=0, \ \qquad \qquad \quad t>0, &x\in\partial D,
    \end{array}\right.\lbl{1.1}\ees
where $\alpha>0$, and $ D $ is a bounded domain in $\mathbb{R}^n$ with smooth boundary
$\partial D$. It was shown (\cite{F}) that for a nonnegative initial condition $u_0\in L^2( D)$
satisfying
   \bes
\int_ D u_0(x)\phi(x)dx>\lambda_1^{\frac{1}{\alpha}},
  \lbl{1.2}\ees
the solution develops  finite time blow-up. Here $\lambda_1$ is the first eigenvalue
of the Laplacian operator $-\Delta$, with zero Dirichlet boundary condition on $ D$, and $\phi$ is the corresponding eigenfunction
normalized so that $\|\phi\|_{L^1( D)}=1$.
 Kaplan \cite{sK} showed that the solution of (\ref{1.1}) will blow up if the initial
datum is large enough.   Fujita \cite{F,F1}  proved that the Cauchy problem (\ref{1.1}), with $ D=\mathbb{R}^n$, has no global positive
nontrivial solutions if $0<\alpha<2/n$, and every solution with arbitrarily small initial datum blows up.
The same is true for $\alpha=2/n$ as     shown by Hayakawa \cite{kH}. When $\alpha>2/n$, solutions
with small initial data tend to zero as time goes to infinity. In this paper,
we will prove   that   noise can lead to finite time blow-up.

For stochastic parabolic equations, the existence of   solutions  has been well studied  \cite{LR,L,T}.
 For instance, for the following equation
   \bes\left\{\begin{array}{lll}
   du=(\Delta u+f(u))dt+\sigma(u)dW_t, \ \ \qquad t>0,&x\in  D,\\[1.5mm]
   u(x,0)=u_0(x)\geq0, \ \ \ &x\in D,\\[1.5mm]
   u(x,t)=0, \qquad \qquad \qquad \qquad \qquad \qquad t>0,  &x\in\partial D,
    \end{array}\right.\lbl{1.3}\ees
Da Prato-Zabczyk \cite{PZ1992} considered  the existence of global solutions   with
additive noise ($\sigma$ is constant). Manthey-Zausinger \cite{MZ} considered (\ref{1.3}),
where $\sigma$ satisfied the global Lipschitz condition. Dozzi and L\'{o}pez-Mimbela
\cite{DL} considered   equation (\ref{1.3}) with
$\sigma(u)=u$. They proved that if $f(u)\geq u^{1+\alpha}$ ($\alpha>0$) and initial data is large enough, the solution will blow up in finite time,
and that if $f(u)\leq u^{1+\beta}$ ($\beta$ is a certain positive constant) and the initial data is
small enough, the solution will exist globally, also see \cite{NX}. A natural question arises: If
$\sigma$ does not satisfy the global Lipschitz condition, what can we say about the solution?
Will it blow up in finite time or exist globally? In a somewhat different case, Mueller \cite{Mu} and, later,
Mueller-Sowers \cite{MuS} investigated the problem of a noise-induced explosion for a
special case of equation (\ref{1.3}), where $f(u)\equiv0,\,\sigma(u)=u^\gamma$ with
$\gamma>0$ and $W(x,t)$ is a space-time white noise. It was shown that the solution will
explode in finite time with positive probability for some $\gamma>3/2$.

\medskip

In the present paper, we shall provide separate
sufficient conditions to ensure that the solutions of (\ref{1.3}) remain positive, or blow up
in finite time. Here blowup means that the solutions will blow up
in finite time in mean $L^\infty$-norm or mean square $L^\infty$-norm; see
Theorems \ref{t3.1}-\ref{t3.3}. Moreover, we will consider a special case, i.e., stochastic Allen-Cahn equation,
whose solution will not blow up in finite time and thus exists globally.

\medskip

This paper is arranged as follows.   After some preliminaries in the next section,
we prove that the solutions of (\ref{2.1}) remain positive under some assumptions in Section 3. Section 4
is concerned with the blow-up phenomenon of solution to (\ref{2.1}) and we will obtain a new result,
which shows that noise can indeed lead to finite time blow-up. In Section 5, we consider
the existence of global solution, with help of a Lyapunov functional,
for a stochastic Allen-Cahn equation, where the noise intensity
$\sigma(u)=u^{1+\beta}$ ($\beta>0$) is not  globally Lipschitz continuous.

\section{Preliminaries}
\setcounter{equation}{0}

To set the stage for our study, we recall Chow's recent works \cite{C2009,C2011} on finite time  blow-up for the following SPDE
      \bes\left\{\begin{array}{lll}
   du=(A u+f(u,x,t))dt+\sigma(u,\nabla u,x,t)dW_t, \ \ \qquad &t>0,\ x\in  D,\\[1.5mm]
   u(x,0)=u_0(x), \ \ \ &\qquad\ \ x\in D,\\[1.5mm]
   u(x,t)=0, \ \ \ \ &t>0,\ x\in\partial D,
    \end{array}\right.\lbl{1.4}\ees
where $A=\sum_{i,j=1}^n\frac{\partial}{\partial_{x_i}}(a_{ij}\frac{\partial}{\partial_{x_j}})$ is
a symmetric, uniformly elliptic operator with smooth coefficients,
$\sigma$ is a given function, and $W(x,t)$ is a Wiener random field defined in a complete probability
space $(\Omega,\mathcal {F},\mathbb{P})$ with a filtration $\mathcal {F}_t$. The Wiener random field has mean
$\mathbb{E}W(x,t)=0$ and its covariance function $q(x,y)$ is defined by
    \bess
\mathbb{E}W(x,t)W(y,s)=(t\wedge s)q(x,y), \ \ \ x,y\in\mathbb{R}^n,
   \eess
where $(t\wedge s)=\min\{t,s\}$ for $0\leq t,s\leq T$. The existence
of strong solutions of (\ref{1.4}) has been studied by many authors  \cite{Cb,PZ}.
To consider positive solutions, they  start with  the unique
solution $u\in C(\bar D\times[0,T])\cap L^2((0,T);H^2)$ for equation (\ref{1.4}).
Under the following conditions
   \begin{quote}
 (P1) There exists a constant $\delta\geq0$ such that
 \bess
 \frac{1}{2} q(x,x)\sigma^2(r,\xi,x,t)-\sum_{i,j=1}^na_{ij}\xi_i\xi_j\leq\delta r^2
    \eess
 for all $r\in\mathbb{R},x\in\bar D,\xi
\in\mathbb{R}^n$ and $t\in[0,T]$;\\
 (P2) The function $f(r,x,t)$ is continuous on $\mathbb{R}\times\bar D\times[0,T]$ and such that
 $f(r,x,t)\geq0$ for $r\leq0$ and $x\in\bar D$, $t\in[0,T]$;  and \\
 (P3) The initial datum $u_0(x)$ on $\bar D$ is positive and continuous,
 \end{quote}
Chow obtained the following result \cite{C2009}.

\begin{prop}\lbl{p1.1}{\rm\cite[Theorem 3.3]{C2009}} Suppose that the conditions
{\rm (P1),(P2)} and {\rm (P3)} hold true. Then the solution of the initial-boundary problem
for the parabolic It\^{o}'s equation {\rm(\ref{1.4})} remains positive, i.e.,
$u(x,t)\geq 0$, a.s. for almost every $x\in D$ and for all $ t\in[0,T]$.
\end{prop}

From (P1), it follows that $\sigma=ku$ ($k$ is a constant) if we only consider the case that $\sigma=\sigma(u)$.
The similar result can be found in \cite{DL,Shiga}. The previous results on existence of global
solution to (\ref{1.4}) require that $\sigma(u)$ satisfies a  global positive Lipschitz condition. A natural
question is what the solution becomes if $\sigma(u)$ does not satisfy the global Lipschitz condition.
We shall study the positive solutions and global solutions of (\ref{1.4}) with $\sigma(u)=u^\gamma$ (for $\gamma>1$)
in sections 3 and 5, respectively.

We consider the eigenvalue problem for the elliptic equation
  \bes\left\{\begin{array}{llll}
-\Delta \phi=\lambda \phi, \ \ \ \ \ \ \  \ {\rm in} \  D,\\
\phi=0, \ \ \qquad\ \qquad  {\rm on}\ \partial D.
   \end{array}\right.\lbl{1.5}\ees
Then, all the eigenvalues are strictly positive, increasing and
the eigenfunction $\phi$ corresponding to the smallest eigenvalue
$\lambda_1$ does not change sign in domain $ D$, as shown in \cite{GT}.
Therefore, we   normalize it in such a way that
   \bess
\phi(x)\geq0,\ \ \ \ \int_ D \phi(x)dx=1.
   \eess
In paper \cite{C2011}, Chow assumed that the following conditions hold, where
$\lambda_1$ is the first eigenvalue of (\ref{1.5}) with $\Delta$ replaced by $A$.

   \begin{quote}
   (N1) There exist a continuous function $F(r)$ and a constant $r_1>0$ such that
   $F$ is positive, convex and strictly increasing for $r\geq r_1$ and satisfies
      \bess
f(r,x,t)\geq F(r)
   \eess
for $r\geq r_1$, $x\in\bar D$, $t\in[0,\infty)$;\\
(N2) There exists a constant $M_1>r_1$ such that $F(r)>\lambda_1r$ for $r\geq M_1$;\\
(N3) The positive initial datum satisfies the condition
   \bess
(\phi,u_0)=\int_ D u_0(x)\phi(x)dx>M_1;
   \eess
(N4) The following condition holds
   \bess
\int_{M_1}^\infty\frac{dr}{F(r)-\lambda_1r}<\infty.
   \eess
    \end{quote}
Alternatively, he imposes the following conditions $S$ on the noise term:
   \begin{quote}
(S1) The correlation function $q(x,y)$ is continuous and positive for $x,y\in\bar D$
such that
   \bess
\int_ D\int_ D q(x,y)v(x)v(y)dxdy\geq q_1\int_ D v^2(x)dx
   \eess
for any positive $v\in H$ and for some $q_1>0$;

(S2) There exist a positive constant $r_2$, continuous functions $\sigma_0(r)$ and $G(r)$ such that
they are both positive, convex and strictly increasing for $r\geq r_2$ and satisfy
  \bess
\sigma(r,x,t)\geq \sigma_0(r)\ \ \ \ {\rm and} \ \ \ \ \sigma_0^2(r)\geq2G(r^2)
   \eess
for $x\in\bar D$, $t\in[0,\infty)$;

(S3) There exists a constant $M_2>r_2$ such that $q_1G(r)>\lambda_1r$ for $r\geq M_2$;

(S4) The positive initial datum satisfies the condition
   \bess
(\phi,u_0)=\int_ D u_0(x)\phi(x)dx>M_2;
   \eess

(S5) The following integral is convergent so that
   \bess
\int_{M_2}^\infty\frac{dr}{q_1G(r)-\lambda_1r}<\infty.
   \eess
    \end{quote}

\begin{prop}\lbl{p1.2} {\rm\cite[Theorem 3.1]{C2011}}
Suppose the initial-boundary value problem {\rm(\ref{1.4})} has a unique local solution
and the conditions {\rm(P1)-(P3)} are satisfied, where $\sigma$ does not depend on
$\nabla u$. In addition, we assume that either
the conditions {\rm(N1)-(N4)} or the alternative conditions {\rm(S1)-(S5)} given above
hold true. Then, for a real number $p>0$, there exists a constant $T_p>0$ such that
   \bes
\lim\limits_{t\rightarrow T_p-}\mathbb{E}\|u \|_p
=\lim\limits_{t\rightarrow T_p-}\mathbb{E}\left(\int_ D|u(x,t)|^pdx\right)^\frac{1}{p}=\infty,
  \lbl{1.6}\ees
or the solution explodes in the mean $L^p$-norm as shown by {\rm(\ref{1.6})},
where $p\geq1$ under conditions $N$, while $p\geq2$ under conditions $S$.
\end{prop}

Looking at the conditions in Propositions \ref{p1.1} and \ref{p1.2}, it is clear that the condition (P1) is very stringent. A noise intensity like $\sigma(u)=u^{1+\beta}$, $\beta>0$,
does not satisfy the condition (P1). But the condition (S5) implies that
$G(r)\geq r^{1+\varepsilon}$, where $\varepsilon$ is a positive constant. Therefore, in order
to prove that noise can lead to blow up, we should delete or change the condition (P1).
Moreover, if we assume that $\sigma=\sigma(u)$, we see that the term
$-\sum_{i,j=1}^na_{ij}(x)\partial_{x_i} u\partial_{x_j} u$ will not play any role in proving that
the solutions are positive.  When the elliptic operator is replaced by the $p$-Laplace operator, the situation is different.

%Inspired by the work \cite{C2009,C2011},

\section{Positive solutions}
\setcounter{equation}{0}

In this section, we will consider the positive solution
to (\ref{2.1}), which will be used to examine the  finite time blow-up phenomenon.

For simplicity, we first consider the following stochastic parabolic It\^{o} equation
    \bes\left\{\begin{array}{lll}
   du=(\Delta u+f(u,x,t))dt+\sigma(u,\nabla u,x,t)dW_t, \ \ \qquad & t>0,\ x\in  D,\\[1.5mm]
   u(x,0)=u_0(x), \ \ \ &\qquad\ \ x\in D,\\[1.5mm]
   u(x,t)=0, \ \ \ \ & t>0,\ x\in\partial D.
    \end{array}\right.\lbl{2.1}\ees
We assume that the covariance function $q(x,y)$
 is bounded, continuous and there is a constant $q_0>0$ such that
    \bess
\sup_{x,y\in D}|q(x,y)|\leq q_0\ \ \ \ {\rm and}\ \ \ \ \int_\mathbb{R}q(x,x)dx<\infty.
   \eess
In addition, we assume that
 \bes\left\{\begin{array}{llll}\medskip
f(u,x,t)\geq a_1u^\beta+a_2 u,\\
\ds\frac{q_0}{2}\sigma^2(u,\nabla u,x,t)-|\nabla u|^2\leq b_1u^{2m}+b_2u^2,
   \end{array}\right.\lbl{2.2}\ees
where $a_2\in\mathbb{R},\,b_i,\beta\geq0$, $(-1)^\beta\in\mathbb{R}$ and
   \bess
&&a_1\left\{\begin{array}{lll}
>0,\ \ \ \ {\rm if}\ (-1)^{\beta}=1,\\[1.5mm]
<0,\ \ \ \ {\rm if}\ (-1)^{\beta}=-1,
   \end{array}\right. \ \ \ \ \  {\rm and} \\[1mm]
&&1\leq m<(1+\beta)/2.
  \eess
As in \cite{C2009,C2011}, let $\eta(r)=r^-$ denote the negative part of $r$ for $r\in\mathbb{R}$.
Set
   \bess
k(r)=\eta^2(r),
   \eess
so that $k(r)=0$ for $r\geq0$ and $k(r)=r^2$ for $r<0$. For
$\varepsilon>0$, let $k_\varepsilon(r)$ be a $C^2$-regularization of $k(r)$ defined by
   \bess
k_\varepsilon(r)=\left\{\begin{array}{lllll}
r^2-\ds\frac{\varepsilon^2}{6}, \ \qquad \ &r<-\varepsilon,\\[1mm]
-\ds\frac{r^3}{\varepsilon}\left(\ds\frac{r}{2\varepsilon}+\frac{4}{3}\right),\ \ &-\varepsilon\leq r<0,\\[1mm]
0, \ \ \ &r\geq0.
   \end{array}\right.\eess
Then one can check that $k_\varepsilon(r)$ has the following properties.
  \begin{lem}\lbl{l2.1}{\rm\cite[Lemma 3.1]{C2009}} The first two derivatives $k'_\varepsilon,\,k''_\varepsilon$
  of $k_\varepsilon$ are continuous and satisfy the conditions: $k'_\varepsilon(r)=0$ for $r\geq0$;
  $k'_\varepsilon\leq0$ and $k''_\varepsilon\geq0$ for any $r\in\mathbb{R}$. Moreover, as $\varepsilon\rightarrow0$,
  we have
    \bess
k_\varepsilon(r)\rightarrow k(r),\ \ k'_\varepsilon(r)\rightarrow-2\eta(r)\ \ \ {\rm and }\
\ k''_\varepsilon(r)\rightarrow2\theta(r),
   \eess
where $\theta(r)=0$ for $r\geq0$, $\theta=1$ for $r<0$, and the convergence is uniform for $r\in\mathbb{R}$.
   \end{lem}
\begin{lem}\lbl{l2.2}{\rm\cite[Lemma 7.6]{GT}} If $u\in W^1( D)$; then $u^+,\,u^-,\,|u|\in W^1( D)$ and
   \bess
\nabla u^-=\left\{\begin{array}{llll}
0,\ \ \ \ \ & u\geq0,\\
\nabla u,\ \ \ \ \ &u<0.
   \end{array}\right.\eess
\end{lem}
With the aid of the above lemmas, we can obtain the following positivity result.

  \begin{theo}\lbl{t2.1} Suppose that {\rm(\ref{2.2})} holds with $1\leq m<(\beta+1)/2$. Then the solution of
initial-boundary value problem {\rm(\ref{2.1})} with nonnegative initial datum remains positive, i.e.,
$u(x,t)\geq 0$, a.s., for almost every $x\in D$ and for all $ t\in[0,T]$.
 \end{theo}

{\bf Proof.} We remark that when $m=1$, Theorem \ref{t2.1} has been proved in \cite{C2009}.
Let $u_t=u(\cdot,t)$ and
   \bess
\Phi_\varepsilon(u_t)=(1,k_\varepsilon(u_t))=\int_ D k_\varepsilon(u(x,t))dx.
   \eess
By It\^{o}'s formula, we have
    \bess
\Phi_\varepsilon(u_t)&=&\Phi_\varepsilon(u_0)+\int_0^t\int_ D k_\varepsilon'(u(x,s))\Delta u(x,s)dxds\\
&&+\int_0^t\int_ D k_\varepsilon'(u(x,s))f(u(x,s),x,s)dxds\\
&&+\int_0^t\int_ D k_\varepsilon'(u(x,s))\sigma(u(x,s),\nabla u(x,s),x,s)dW(x,s)dx\\
&&+\frac{1}{2}\int_0^t\int_ D k_\varepsilon''(u(x,s))q(x,x)\sigma^2(u(x,s),\nabla u(x,s),x,s)dxds\\
&=&\Phi_\varepsilon(u_0)+\int_0^t\int_ D k_\varepsilon''(u(x,s))\left(\frac{1}{2}q(x,x)\sigma^2(u(x,s),\nabla u(x,s),x,s)-|\nabla u|^2\right)dxds\\
&&+\int_0^t\int_ D k_\varepsilon'(u(x,s))f(u(x,s),x,s)dxds\\
&&+\int_0^t\int_ D k_\varepsilon'(u(x,s))\sigma(u(x,s),\nabla u(x,s),x,s)dW(x,s)dx.
   \eess
Taking expectation over the above equality, we get
    \bess
\mathbb{E}\Phi_\varepsilon(u_t)
&=&\Phi_\varepsilon(u_0)+\mathbb{E}\int_0^t\int_ D k_\varepsilon''(u(x,s))\\
&&\times\left(\frac{1}{2}q(x,x)\sigma^2(u(x,s),\nabla u(x,s),x,s)-|\nabla u|^2\right)dxds\\
&&+\mathbb{E}\int_0^t\int_ D k_\varepsilon'(u(x,s))f(u(x,s),x,s)dxds.
  \eess
Note that $\lim\limits_{\varepsilon\rightarrow0}\mathbb{E}\Phi_\varepsilon(u_t)=\mathbb{E}\|\eta(u_t)\|^2$, by
taking the limits termwise as $\varepsilon\rightarrow0$ and using Lemma \ref{l2.1}, we have
    \bes
\mathbb{E}\|\eta(u_t)\|^2&=&\mathbb{E}\|\eta(u_0)\|^2
+2\mathbb{E}\int_0^t\int_ D \nm\\
&&\left(\frac{1}{2}q(x,x)\sigma^2(u^-(x,s),\nabla u^-(x,s),x,s)-|\nabla u^-|^2\right)dxds\nm\\
&&-2\mathbb{E}\int_0^t\int_ D \eta(u(x,s))f(u(x,s),x,s)dxds,
  \lbl{2.3}\ees
where $\|\cdot\|$ denotes the norm of $L^2( D)$. We remark that $\nabla u^-$ exists by Lemma \ref{l2.2}.
Moreover, it follows from Lemma \ref{l2.2} that (\ref{2.3}) is well defined.
The authors  of \cite{LR} proved that $u\in L^p( D)$ if
$u_0\in L^p( D)$, where $p\geq1$ and $u$ is the solution of (\ref{1.3}).
We also remark that in \cite{LR} they assumed $\sigma$ satisfied the linear
growth. Taniguchi \cite{T} obtained the well-posedness of (\ref{1.3}) under the
condition that the nonlinear terms $f$ and $\sigma$ satisfy the local Lipschitz condition.
 Chow \cite{Cb} obtained the well
posedness of (\ref{2.1}), where the nonlinear terms $f$ and $\sigma$ satisfy the global Lipschitz
condition (see p.74-84 in \cite{Cb}).
One can use the method of \cite{T} to obtain the well-posedness of (\ref{2.1}),
where $f$ and $\sigma$ satisfy the local Lipschitz condition. By using (\ref{2.2}) and $\eta(u)\geq0$,
we obtain
   \bes
\mathbb{E}\|\eta(u_t)\|^2
&\leq&\mathbb{E}\|\eta(u_0)\|^2-2\mathbb{E}\int_0^t\int_ D \eta(u(x,s))(a_1u^\beta(x,s)+a_2u)dxds\nm\\
&&+2\mathbb{E}\int_0^t\int_ D
\left(\frac{1}{2}q_0\sigma^2(u^-(x,s),\nabla u^-(x,s),x,s)-|\nabla u^-|^2\right)dxds\nm\\
&\leq&\mathbb{E}\|\eta(u_0)\|^2-2\mathbb{E}\int_0^t\int_ D (|a_1|(u^-)^{\beta+1}(x,s)-a_2(u^-)^2(x,s))dxds\nm\\
&&+2\mathbb{E}\int_0^t\int_ D
\left(b_1(u^-)^{2m}(x,s)+b_2(u^-)^2(x,s)\right)dxds,
  \lbl{j.1}\ees
where we have used the condition on $a_1$, that is, $(-1)^\beta a_1=|a_1|$.

Now, we   use $\epsilon$-Young's
inequality and the following interpolation
inequality of $L^p$ to deal with the last two terms of (\ref{j.1}),
  \bes
\|u\|_{L^r}\leq\|u\|^\theta_{L^p}\|u\|^{1-\theta}_{L^q},
   \lbl{2.4}\ees
where $0<\theta<1$ and
   \bes
\frac{1}{r}=\frac{\theta}{p}+\frac{1-\theta}{q},\ \ \ 0<p<q\leq\infty.
   \lbl{2.5}\ees
For simplicity, we write
$u$ instead of $u^-$. Notice that $1<m<\frac{\beta+1}{2}$ and $0\leq q(x,x)\leq q_0$,
by using (\ref{2.4}), we have
   \bes
2b_1\int_ D u^{2m}(x,t)dx&=&2b_1\|u\|^{2m}_{L^{2m}}\nm\\
&\leq&C\|u\|_{L^2}^{2m\theta}\|u\|_{L^{\beta+1}}^{2m(1-\theta)}\nm\\
&\leq&\epsilon\|u\|_{L^{\beta+1}}^{2m(1-\theta)\frac{1}{1-m\theta}}+C(\epsilon)\|u\|_{L^2}^{2}\nm\\
&=&\epsilon\|u\|_{L^{\beta+1}}^{\beta+1}+C(\epsilon)\|u\|_{L^2}^{2},
  \lbl{2.6}\ees
where $\theta=\frac{\beta+1-2m}{m(\beta-1)}$ satisfying (\ref{2.5}). Substituting (\ref{2.6}) into (\ref{j.1})
and using the fact $\|\eta(u_0)\|=0$ because $u_0$ is non-negative,
we get
    \bess
\mathbb{E}\|\eta(u_t)\|^2\leq2\int_0^t\mathbb{E}(\epsilon\|u^-_s\|_{L^{\beta+1}}^{\beta+1}+(C(\epsilon)+b_2+a_2)\|u^-_s\|_{L^2}^{2})ds
-2|a_1|\int_0^t\mathbb{E}\|u^-_s\|_{L^{\beta+1}}^{\beta+1}ds.
  \eess
Let $0<\epsilon<|a_1|$. We observe that
    \bess
\mathbb{E}\|\eta(u_t)\|^2\leq C\int_0^t\mathbb{E}\|\eta(u_s)\|_{L^2}^{2}ds,
  \eess
which, by means of Gronwall's inequality, implies that
  \bess
\mathbb{E}\|\eta(u_t)\|^2=0, \ \ \ \ \forall t\in[0,T].
   \eess
It follows that $\eta(u_t)=u^-(x,t)=0$ a.s. for a.e. $x\in D$, $\forall t\in[0,T]$.
This completes the proof. $\Box$

%Secondly, we consider the other case $0<m\leq1$.
%  \begin{theo}\lbl{t2.2} Suppose that $0<m\leq1$ and the following inequality hold:
%   \bes
%(\eta(u),f(u))\geq-\delta\|\eta(u)\|_{L^2}^2,
%    \lbl{2.7}\ees
%where $\delta>0$. Then the solution of
%initial-boundary value problem {\rm(\ref{2.1})} with nonnegative initial data remains positive so that
%$u(x,t)\geq0$, a.s. for almost every $x\in D$, $\forall t\in[0,T]$.
% \end{theo}

%{\bf Proof.}

\begin{remark}\lbl{r2.1} $1.$ Comparing Theorem {\rm\ref{t2.1}} with Proposition {\rm \ref{p1.1}},
it is clear  that our assumption is weaker. For example,   $f(u)=u(1-u^2)$,
  will not satisfy the condition {\rm(P2)}, but it is covered in our theorem. By using a similar method,
one can deal with the nonlinearity term depending on the $x$ and $t$. In this section, we only consider
the case that $m>1$ and it is possible to use the similar method to deal with the case that $0<m<1$.

$2.$ Obviously, if $f(u)\equiv0$, Theorem {\rm\ref{t2.1}} will fail, that is, we can not obtain
the positivity of solutions to {\rm(\ref{2.1})} with $f(u)\equiv0$, even for one dimension. Because
in the proof Theorem {\rm\ref{t2.1}}, we use $f(u)$ to control $\sigma(u)$. But we can
obtain the positivity of solutions of {\rm(\ref{2.1})} with $f(u)\equiv0$ under the condition that
$\Delta$ is replaced by the $p-$Laplacian $\Delta_p$.

$3.$ From the proof of Theorem {\rm\ref{t2.1}}, we know that the term $-|\nabla u|^2$ is
  a good term, as we can use
it to control the stochastic term. But when we use the embedding theorem and interpolation
inequality, we find that the term $\|u\|_{L^m}^m$ would be changed to $\|u\|_{L^2}^\nu$, where
$\nu>2$. Due to the convexity of the function $x^\nu$, we can not get the desired result.
However, if we change $-|\nabla u|^2$ to  $-|\nabla u|^p$, $p>n\geq2$, we can show that the solution
is positive; see     Theorem {\rm\ref{t2.2}} below.
\end{remark}

\medskip

It follows from Theorem \ref{t2.1} that the value of $m$ depends on the nonlinear term $f$. The
following result shows that the value of $m$ may not depend on the nonlinearity term $f$.
Now, we consider the following It\^{o} parabolic equation with the $p-$Laplacian operator
    \bes\left\{\begin{array}{lll}
   du=(\Delta_p u+f(u,x,t))dt+g(u,x,t)dW_t, \ \ \qquad & t>0,\ x\in  D,\\[1.5mm]
   u(x,0)=u_0(x), \ \ \ &\qquad\ \ x\in D,\\[1.5mm]
   u(x,t)=0, \ \ \ \ & t>0,\ x\in\partial D,
    \end{array}\right.\lbl{2.7}\ees
where $\Delta_pu=div(|\nabla u|^{p-2}\nabla u)$. We assume that there exist positive constants
$\alpha$ and $\beta$ such that
   \bes
&&(\eta(v),f(v,x,t))\geq-\alpha\|\eta(v)\|^2,\nm\\
&&g^2(u,x,t)\leq 2\beta u^{2m},\ \ \ m>1.
   \lbl{2.8}\ees

\begin{theo}\lbl{t2.2} Assume that $p>\max\{2m,n\}$ and {\rm(\ref{2.8})} holds. Then the solution of
initial-boundary value problem {\rm(\ref{2.7})} with nonnegative initial datum remains positive:
$u(x,t)\geq0$, a.s., for almost every $x\in D$ and all $ t\in[0,T]$.
 \end{theo}

{\bf Proof.} Similar to the proof of Theorem \ref{t2.1}, we   get
    \bes
\mathbb{E}\|\eta(u_t)\|^2&=&E\|\eta(u_0)\|^2
+2\mathbb{E}\int_0^t\int_ D \left(\frac{1}{2}q(x,x)g(u^-,x,t)-|\nabla u^-|^p\right)dxds\nm\\
&&-2\mathbb{E}(\eta(v),f(u,x,t))\nm\\
&\leq&\mathbb{E}\|\eta(u_0)\|^2
+2\mathbb{E}\int_0^t\int_ D \left(\beta q_0(u^-)^{2m}-|\nabla u^-|^p\right)dxds\nm\\
&&+2\alpha E\|u^-\|^2.
  \lbl{2.9}\ees
It follows from Lemma \ref{l2.2} that the above inequality is well defined. For simplicity,
we write $u$ instead of $u^-$. By the
Sobolev embedding inequality and  for $p>n$,
   \bes
\|u\|_{L^\infty}\leq C\|u\|_{W^{1,p}},
   \lbl{2.10}\ees
which implies that
    \bes
\|u\|_{L^{2m}}^{2m}=\int_ D u^{2m}(x,t)dx&\leq&\|u\|_{L^\infty}^{2m-2}\int_ D u^2(x,t)dx\nm\\
&=&\|u\|_{L^\infty}^{2m-2}\|u\|_{L^2}^{2}\nm\\
&\leq&\|u\|_{L^\infty}^{2m-2+\gamma}\|u\|_{L^2}^{2-\gamma}\nm\\
&\leq&C\|u\|_{W^{1,p}}^{2m-2+\gamma}\|u\|_{L^2}^{2-\gamma}\nm\\
&\leq&C(\varepsilon)\|u\|_{L^2}^{2}+\varepsilon\|u\|_{W^{1,p}}^{(2m-2+\gamma)\cdot\frac{2}{\gamma}},
  \lbl{2.11}\ees
where $\varepsilon>0$ and we have used the $\epsilon$-Young's inequality. Noting that $p>\max\{2m,n\}$,
there exists a constant $\gamma\in(0,2)$ such that
   \bess
(2m-2+\gamma)\cdot\frac{2}{\gamma}=p.
   \eess
Letting $kq_0\varepsilon\leq1$ and submitting (\ref{2.11}) into (\ref{2.9}), we get
    \bess
\mathbb{E}\|\eta(u_t)\|^2\leq C\int_0^t\mathbb{E}\|\eta(u_s)\|_{L^2}^{2}ds,
  \eess
which, by means of Gronwall's inequality, implies that
  \bess
\mathbb{E}\|\eta(u_t)\|^2=0, \ \ \ \ \forall t\in[0,T].
   \eess
It follows that $\eta(u_t)=u^-(x,t)=0$ a.s. for a.e. $x\in D$, $\forall t\in[0,T]$.
This completes the proof.
$\Box$

We remark that the value of $m$ in Theorems \ref{t2.1} and \ref{t2.2} either depends on
the nonlinear term $f$ or the operator $\Delta_p$.
% The reason is that we didn't choose suitable test function.
In the followings, we will select a new test function $\beta_\varepsilon(r)$,
instead of $k_\varepsilon(r)$.

Define new functions
    \bes
&&\beta_\varepsilon(r)=\int_r^\infty\rho_\varepsilon(s)ds,\ \ \
\rho_\varepsilon(r)=\int_{r+\varepsilon}^\infty J_\varepsilon(s)ds,\ \ \ r\in\mathbb{R},\nm\\
&&J_\varepsilon(x)=\varepsilon^{-n}J\left(\frac{|x|}{\varepsilon}\right),\ \ \
J(x)=\left\{\begin{array}{llll}\medskip
C\exp\left(\frac{1}{|x|^2-1}\right),\ \ \ &|x|<1,\\
0, \ \ \ &|x|\geq1.
  \end{array}\right.\lbl{2.13}\ees
Then by direct verification, we have the following result.

\begin{lem}\lbl{l2.3} The above constructed functions $\rho_\varepsilon,\beta_\varepsilon$ are in $C^\infty(\mathbb{R})$ and
have the following properties:
$\rho_\varepsilon$ is a non-increasing function and
   \bess
\beta_\varepsilon'(r)=-\rho_\varepsilon(r)=\left\{\begin{array}{llll}\medskip
0,\ \ \ &r\geq0,\\
1,\ \ \ &r\leq-2\varepsilon.
 \end{array}\right.\eess
 Additionally,  $\beta_\varepsilon$ is convex and
    \bess
\beta_\varepsilon(r)=\left\{\begin{array}{llll}\medskip
0,\ \  \ \ \ &r\geq0,\\
-2\varepsilon-r+\varepsilon \hat C, \ \  &r\leq-2\varepsilon,
   \end{array}\right.\eess
where $\hat C=\int^0_{-2}\int_{t+1}^1J(s)dsdt<2$. Furthermore,
   \bess
0\leq\beta_\varepsilon''(r)=J_\varepsilon(r+\varepsilon)\leq \varepsilon^{-n}C, \ \ -2\varepsilon\leq r\leq0,
   \eess
which implies that
    \bess
-2^nC\leq r^n\beta_\varepsilon''(r)\leq0 \ \ \ \ & {\rm for}\ -2\varepsilon\leq r\leq0,\ {\rm and}\ n\ {\rm is}\ {\rm odd};\\
0\leq r^n\beta_\varepsilon''(r)\leq2^nC\ \ \ \ &{\rm for}\ -2\varepsilon\leq r\leq0,\ {\rm and}\ n\ {\rm is}\ {\rm even}.
   \eess
   \end{lem}

Now, we consider the following stochastic parabolic It\^{o} equation
    \bes\left\{\begin{array}{lll}
   du=(\Delta u+f(u,x,t))dt+g(u,x,t)dW_t, \ \ \qquad & t>0,\ x\in  D,\\[1.5mm]
   u(x,0)=u_0(x), \ \ \ &\qquad\ \ x\in D,\\[1.5mm]
   u(x,t)=0, \ \ \ \ & t>0,\ x\in\partial D.
    \end{array}\right.\lbl{2.14}\ees

\begin{theo}\lbl{t2.3} Assume that (i) the function $f(r,x,t)$ is continuous on
  $\mathbb{R}\times\bar D\times[0,T]$;  (ii) $f(r,x,t)\geq0$ for $r\leq0$, $x\in\bar D$ and $t\in[0,T]$; and (iii)  $ g^2(u,x,t)\leq ku^{2m}$,
where $k>0$, $2m>n$ and $(-1)^{2m-n}\in\mathbb{R}$. Then the solution of
initial-boundary value problem {\rm(\ref{2.14})} with nonnegative initial datum remains positive:
$u(x,t)\geq0$, a.s. for almost every $x\in D$ and for all $  t\in[0,T]$.
\end{theo}

{\bf Proof.}  Define
   \bess
\Phi_\varepsilon(u_t)=(1,\beta_\varepsilon(u_t))=\int_ D \beta_\varepsilon(u(x,t))dx.
   \eess
By It\^{o}'s formula, we have
    \bess
\Phi_\varepsilon(u_t)&=&\Phi_\varepsilon(u_0)+\int_0^t\int_ D \beta_\varepsilon'(u(x,s))\Delta u(x,s)dxds\\
&&+\int_0^t\int_ D \beta_\varepsilon'(u(x,s))f(u(x,s),x,s)dxds\\
&&+\int_0^t\int_ D \beta_\varepsilon'(u(x,s))g(u(x,s),x,s)dW(x,s)dx\\
&&+\frac{1}{2}\int_0^t\int_ D \beta_\varepsilon''(u(x,s))q(x,x)g^2(u(x,s),x,t)dxds\\
&=&\Phi_\varepsilon(u_0)+\int_0^t\int_ D \beta_\varepsilon''(u(x,s))\left(\frac{1}{2}q(x,x)g^2(u(x,s),x,s)-|\nabla u|^2\right)dxds\\
&&+\int_0^t\int_ D \beta_\varepsilon'(u(x,s))f(u(x,s),x,s)dxds\\
&&+\int_0^t\int_ D \beta_\varepsilon'(u(x,s))g(u(x,s),x,s)dW(x,s)dx.
   \eess
Taking expectation over the above equality and using Lemma \ref{l2.3}, we get
    \bess
\mathbb{E}\Phi_\varepsilon(u_t)
&=&\mathbb{E}\Phi_\varepsilon(u_0)+\mathbb{E}\int_0^t\int_ D \beta_\varepsilon''(u(x,s))\\
&&\times\left(\frac{1}{2}q(x,x)g^2(u(x,s),x,s)-|\nabla u|^2\right)dxds\\
&&+\mathbb{E}\int_0^t\int_ D \beta_\varepsilon'(u(x,s))f(u(x,s),x,s)dxds\\
&\leq&\mathbb{E}\Phi_\varepsilon(u_0)+\frac{k}{2}\mathbb{E}\int_0^t\int_ D \beta_\varepsilon''(u(x,s))
q(x,x)u(x,s)^{2m}dxds\\
&&+\mathbb{E}\int_0^t\int_ D \beta_\varepsilon'(u(x,s))f(u(x,s),x,s)dxds.
  \eess
Here and after, we denote $\|\cdot\|_{L^1}$ by $\|\cdot\|_1$. Let $\eta(u)=u^-$
denote the negative part of $u$ for $u\in\mathbb{R}$.
Then we have $\lim\limits_{\varepsilon\rightarrow0}\mathbb{E}\Phi_\varepsilon(u_t)=\mathbb{E}\|\eta(u_t)\|_1$.
It follows from Lemma \ref{l2.3} that
    \bess
0\geq u^{2m}\beta''_\varepsilon(u)\geq\left\{\begin{array}{llll}\medskip
0,\ \ \ \ & u\geq0\ {\rm or}\ u\leq-2\varepsilon,\\
-2Cu^{2m-n},\ \ \ \ &-2\varepsilon\leq u\leq0,\ {\rm and}\ u^{2m-1}\geq0,
   \end{array}\right.
   \eess
or
    \bess
0\leq u^{2m}\beta''_\varepsilon(u)\leq\left\{\begin{array}{llll}\medskip
0,\ \ \ \ & u\geq0\ {\rm or}\ u\leq-2\varepsilon,\\
-2Cu^{2m-n},\ \ \ \ &-2\varepsilon\leq u\leq0,\ {\rm and}\ u^{2m-1}\leq0
   \end{array}\right.
   \eess
which implies that
$\lim\limits_{\varepsilon\rightarrow0}u^{2m}\beta_\varepsilon''(u)=0$ provided that $2m>1$.
By taking the limits termwise as $\varepsilon\rightarrow0$ and using Lemma \ref{l2.3}, we get
    \bes
\mathbb{E}\|\eta(u_t)\|_1&\leq&\mathbb{E}\|\eta(u_0)\|_1
-\mathbb{E}\int_0^t\int_ D \eta(u(x,s))f(u(x,s),x,s)dxds\nm\\
&\leq&0,
  \lbl{2.15}\ees
which implies that $u^-=0$ a.s. for a.e. $x\in D$, $\forall t\in[0,T]$.
This completes the proof.
$\Box$

%\begin{remark}\lbl{r2.2} The assumption ``$(-1)^{2m-n}$ makes sense",
%in the Theorem {\rm\ref{t2.3}}, means that $(-1)^{2m-n}$ equals to a real constant.
%There exists some number $m$ and $n$ such that $(-1)^{2m-n}$ does not exist in real number
%domain. For example, let $m=\frac{7}{4}$ and $n=1$, then $(-1)^{2m-1}=(-1)^{5/2}$ will not
%exist in real number domain. On the other hand, $(-1)^{2m-n}$ makes sense for all
%$m\in \mathbb{N}$.

%One can use the same test function to improve the result of Theorem {\rm\ref{t2.2}}.
 %\end{remark}

\section{Blow-up Phenomenon}
\setcounter{equation}{0}

In this section, we shall consider the solutions of (\ref{2.1}) which blow up in
finite time. We first show that a similar result to \cite{F} holds for (\ref{2.1}),
and then we examine how noise induces blow-up
in finite time in the mean $L^\infty$-norm. We divide this section into three subsections.

\subsection{First result on blow-up}

In this subsection, we shall prove that the solution of stochastic parabolic
It\^{o} equation will blow up in finite time if the solution of corresponding
deterministic equation blows up in finite time. Specifically, there exists a
finite time $T^*$ such that $\lim\limits_{t\rightarrow T^*-0}\mathbb{E}\sup_{x\in D}u(x,t)=\infty$,
where $u(x,t)$ is a positive solution of the stochastic parabolic
It\^{o} equation (\ref{2.1}).
We remark that when $\sigma\equiv0$, then (\ref{2.1}) becomes the deterministic
parabolic equation. Indeed,  Fujita \cite{F} presented an   existence and non-existence theorem for
global solution of (\ref{2.1}) with $\sigma\equiv0$. The following result is similar to
that in \cite{F}.
  \begin{theo}\lbl{t3.1} Suppose the initial-boundary value problem {\rm(\ref{2.1})} has
  a unique local solution. Assume that all the assumptions in Theorem
  {\rm\ref{t2.1}} hold, where $a_1>0$. In addition, if $\lambda_1\geq a_2$, we assume that
   \bes
\int_ D u_0(x)\phi(x)dx>[a_1^{-1}(\lambda_1-a_2)]^{\frac{1}{\beta}},
  \lbl{3.1}\ees
and if $\lambda_1<a_2$, we assume that $u_0(x)\geq0$ and $u_0(x)\not\equiv0$, where
$\lambda_1$ is the smallest eigenvalue of the operator $\Delta$ on $D$ and $\phi$
is the corresponding eigenfunction.
Then, there exists a constant $T^*>0$ such that
     \bes
\lim\limits_{t\rightarrow T^*-}\mathbb{E}\|u_t\|_{L^\infty}
=\lim\limits_{t\rightarrow T^*-}\mathbb{E}\sup_{x\in D}|u(x,t)|=\infty.
    \lbl{3.2}\ees
\end{theo}

{\bf Proof.} It follows from Theorem \ref{t2.1} that
(\ref{2.1}) has a unique positive solution. We will prove the theorem by contradiction.
Suppose (\ref{3.2}) is false. Then there exists a global positive solution $u$ such
that for any $T>0$
   \bes
\sup_{0\leq t\leq T}\mathbb{E}\sup_{x\in D}|u(x,t)|<\infty,
  \lbl{3.3}\ees
which implies that
   \bes
\sup_{0\leq t\leq T}\mathbb{E}\int_ D u(x,t)\phi(x)dx\leq\sup_{0\leq t\leq T}\mathbb{E}\sup_{x\in D}|u(x,t)|<\infty,
   \lbl{3.4}\ees
where $\phi$ is defined as in (\ref{1.5}) and satisfies $\int_ D\phi(x)dx=1$.
Define
   \bess
\hat u(t):=\int_ D u(x,t)\phi(x)dx.
   \eess
Then  we have
   \bes
\hat u(t)&=&(u_0,\phi)+\int_0^t\int_ D \Delta u(x,s)\phi(x)dxds+
\int_0^t\int_ D f(u,x,s)\phi(x)dxds\nm\\
&&+\int_0^t\int_ D \sigma(u,\nabla u,x,s)\phi(x)dW(x,s)dx\nm\\
&=&(u_0,\phi)-\lambda_1\int_0^t\int_ D u(x,s)\phi(x)dxds+
\int_0^t\int_ D f(u,x,s)\phi(x)dxds\nm\\
&&+\int_0^t\int_ D \sigma(u,\nabla u,x,s)\phi(x)dW(x,s)dx.
   \lbl{3.5}\ees
Taking the expectation over (\ref{3.5}) and appealing to Fubini's theorem, we obtain
   \bess
\mathbb{E}\hat u(t)=(u_0,\phi)-\lambda_1\int_0^t\mathbb{E}\hat u(s)ds+
\int_0^t\mathbb{E}\int_ D f(u,x,s)\phi(x)dxds,
   \eess
or, in the differential form, for $\xi(t)=\mathbb{E}\hat u(t)$,
   \bes\left\{\begin{array}{llll}
\ds\frac{d\xi(t)}{dt}=-\lambda_1\xi(t)+\mathbb{E}\ds\int_ D f(u,x,t)\phi(x)dx\\[1.5mm]
\xi(0)=\xi_0,
    \end{array}\right.\lbl{3.6}\ees
where $\xi_0=(u_0,\phi)$. By Jensen's inequality, (\ref{3.6})
yields that
     \bes\left\{\begin{array}{llll}
\ds\frac{d\xi(t)}{dt}\geq-\lambda_1\xi(t)+a_1\xi^{\beta}(t)+a_2\xi(t)\\[1.5mm]
\xi(0)=\xi_0.
    \end{array}\right.\lbl{3.7}\ees
For $\xi_0>(a_1\lambda_1)^{\frac{1}{\beta-1}}$, this implies that $a_1\xi^{\beta}(t)-(\lambda_1-a_2)\xi(t)>0$
and $\xi(t)>\xi_0$ for $t>0$. An integration of equation (\ref{3.7}) gives that
   \bess
T\leq\int_{\xi_0}^{\xi(T)}\frac{dr}{a_1r^{\beta}-(\lambda_1-a_2)r}
\leq\int_{\xi_0}^\infty\frac{dr}{a_1r^{\beta}-(\lambda_1-a_2)r}<\infty,
    \eess
which implies $\xi(t)$ must blow up at a time $T^*\leq\int_{\xi_0}^\infty\frac{dr}{a_1r^{\beta}-(\lambda_1-a_2)r}$.
Hence this is a contradiction to (\ref{3.4}). This completes the proof. $\Box$

It is remarked that Proposition \ref{p1.2} covers a part of the above result.
The following example shows that Theorem \ref{t3.1} generalizes Proposition \ref{p1.2}.

{\bf Example} Consider the following stochastic parabolic It\^{o} equation
   \bes\left\{\begin{array}{lll}
   du=(\Delta u+u^{2})dt+ku^{1+\frac{1}{3}}dW_t, \ \ \qquad &t>0,\ x\in  D,\\[1.5mm]
   u(x,0)=u_0(x), \ \ \ &\qquad\ \ x\in D,\\[1.5mm]
   u(x,t)=0, \ \ \ \ &t>0,\ x\in\partial D,
    \end{array}\right.\lbl{3.8}\ees
where $k\in\mathbb{R}$ and $ D$ is defined as in (\ref{1.1}). Fujita \cite{F}
obtained that the solution of (\ref{3.8}) with $k=0$ and $u_0\geq0$ will blow up
in finite time. By Theorem \ref{t2.1},
we know that the solution of (\ref{3.8}) remains positive if $u_0\geq0$. It follows from
Theorem \ref{t3.1} that the solution of (\ref{3.8}) will blow up in finite time
under the same assumptions as in \cite{F}. We
also remark that Proposition \ref{p1.2} is not suitable for (\ref{3.8}).

We also have  the following remarks.

  \begin{remark}\lbl{r3.1}
$1.$ From the proof of Theorem {\rm\ref{t3.1}}, we conclude that
the stochastic term does not play a role because the first moment of white noise
is zero.
%In other words, the first moment does not touch the white noise. Clearly,
White noise can not
prevent the blow-up of the solution. If we want to study whether the noise can prevent singularities
$($see {\rm\cite{FF}}$)$, perhaps we should consider the colored noise or complex noise.

%$2.$ In the proof of Theorem {\rm\ref{t3.1}}, the positivity of solutions
%can assure that $\hat u(t)>0$ and so it is important. Combining Proposition
%{\rm\ref{p1.2}} and the fact that
%  \bess
%\lim\limits_{p\rightarrow\infty}\|u\|_{L^p}=\|u\|_{L^\infty},
%   \eess
%one can obtain a similar result to Theorem {\rm\ref{t3.1}} under the assumptions of Proposition
%{\rm\ref{p1.2}}.

$2.$ In {\rm\cite{DL}}, the authors obtained a similar result to Theorem {\rm\ref{t3.1}}. They
assumed that the nonlinearity $f(u)\geq u^{1+\beta}$ $(\beta>0)$ and $\sigma(u,\nabla u,x,t)=u$.

$3.$ If $a_1=1$ and $a_2=0$, then the condition {\rm(\ref{3.1})} becomes {\rm(\ref{1.2})}. That is,
under the same conditions on initial data, the solutions of {\rm(\ref{1.1})} and {\rm(\ref{2.1})}
will blow up in finite time. Thus we can say we obtain a similar result to {\rm\cite{F}}.
  \end{remark}

\subsection{Second result on blow-up}

In this subsection, we consider the issue about how noise may induce finite time blow-up   of the solution of stochastic
partial differential equations.
%From the Remark \ref{r3.1}, it seems that white noise does not play important role. So in this subsection, we
%consider the second moment, which will be used in the proof of the following Theorem \ref{t3.2}.

Consider the following stochastic parabolic It\^{o} equation
   \bes\left\{\begin{array}{lll}
   du=(\Delta u+|u|^{1+\alpha})dt+bu^mdW_t, \ \ \qquad &t>0,\ x\in  D,\\[1.5mm]
   u(x,0)=u_0(x), \ \ \ &\qquad\ \ x\in D,\\[1.5mm]
   u(x,t)=0, \ \ \ \ &t>0,\ x\in\partial D,
    \end{array}\right.\lbl{3.9}\ees
where $b\in\mathbb{R}$, $\alpha>0$ and $1\leq m<1+\frac{\alpha}{2}$. When $m=1$, Dozzi and L\'{o}pez-Mimbela
\cite{DL} obtained the global solution of (\ref{3.9}) if the initial data and the noise are
small enough (see Theorem 5 in \cite{DL}), which is similar to the deterministic case \cite{F}.
It is known that when $b=0$ and the nonnegative initial datum is small
enough, (\ref{3.9}) has a unique global solution \cite{F}. In this subsection, we
will show that noise can induce blow-up.

\begin{theo}\lbl{t3.2} Assume that $u_0$ is a nonnegative continuous function and
   \bes
\inf_{x,y\in D}q(x,y)\geq q_1,\ \ r^{1+\frac{\alpha}{2}}+\frac{b^2q_1}{2}r^m-\lambda_1r>0,\  \ r=\left(\int_ D u_0(x)\phi(x)dx\right)^2,
   \lbl{3.10}\ees
where $\lambda_1$ is defined as in {\rm(\ref{1.5})} and $q(x,y)$ is the
correlation function. Then the solution of {\rm(\ref{3.9})} will blow up in finite time in
$L^2(\Omega)\times L^\infty(D)$, that is,
there exists a constant $T^*>0$ such that
     \bes
\lim\limits_{t\rightarrow T^*-0}\left(\mathbb{E}\|u\|^2_{L^\infty(D)}\right)^{\frac{1}{2}}
=\lim\limits_{t\rightarrow T^*-0}\left(\mathbb{E}\sup_{x\in D}u(x,t)^2\right)^{\frac{1}{2}}=\infty.
    \lbl{3.11}\ees
\end{theo}

{\bf Proof.} By \cite{Cb,LR,PZ}, we know that (\ref{3.9}) has a unique local
solution. It follows from Theorem \ref{t2.1} that the solution of (\ref{3.9})
remains positive. Since
  \bess
\mathbb{E}\hat u^2(t)\leq\mathbb{E}\sup_{x\in D}u(x,t)^2,
   \eess
it suffices to show that
$\mathbb{E}\hat u^2(t)$ blows up in finite time, where $\hat u(t)=(u,\phi)$.

By applying It\^{o}'s formula to $\hat u^2(t)$ and making use of (\ref{1.5}), we get
   \bes
\hat u^2(t)&=&(u_0,\phi)^2-2\lambda_1\int_0^t\hat u^2(s)ds+2
\int_0^t\int_ D\hat u(s) u^{1+\alpha}(x,s)\phi(x)dxds\nm\\
&&+2b\int_0^t\int_ D \hat u(s)u^m(x,s)\phi(x)dW(x,s)ds\nm\\
&&+b^2\int_0^t\int_ D\int_ D q(x,y) u^m(x,s)\phi(x)u^m(y,s)\phi(y)dxdyds
   \lbl{3.12}\ees
Let $\eta(t)=\mathbb{E}\hat u^2(t)$. By taking an expectation over (\ref{3.12}), we obtain
  \bes
\eta(t)&=&(u_0,\phi)^2-2\lambda_1\int_0^t\eta(s)ds
+2\mathbb{E}\int_0^t\int_ D\hat u(s) u^{1+\alpha}(x,s)\phi(x)dxds\nm\\
&&+b^2\int_0^t\mathbb{E}\int_ D\int_ D q(x,y) u^m(x,s)\phi(x)u^m(y,s)\phi(y)dxdyds,
   \lbl{3.13}\ees
or, in the differential form
   \bes\left\{\begin{array}{lllll}
\ds\frac{d\eta(t)}{dt}=-2\lambda_1\eta(t)+2\mathbb{E}\hat u(t)\ds\int_ D u^{1+\alpha}(x,t)\phi(x)dx\\[2.5mm]
\quad\quad\quad\ \ +b^2\mathbb{E}\ds\int_ D\int_ D q(x,y) u^m(x,t)\phi(x)u^m(y,t)\phi(y)dxdy\\[2mm]
\eta(0)=\eta_0=(u_0,\phi)^2.
    \end{array}\right.\lbl{3.14}\ees
By Jensen's inequality, (\ref{3.14})
yields
     \bes\left\{\begin{array}{llll}
\ds\frac{d\eta(t)}{dt}\geq-2\lambda_1\eta(t)+2\eta^{1+\frac{\alpha}{2}}(t)+q_1b^2\eta^m(t)\\[1.5mm]
\eta(0)=\eta_0.
    \end{array}\right.\lbl{3.15}\ees
This implies that, for $\eta_0^{1+\frac{\alpha}{2}}+\frac{b^2q_1}{2}\eta_0^m-\lambda_1\eta_0>0$,
we have $\eta^{1+\frac{\alpha}{2}}(t)+\frac{b^2q_1}{2}\eta(t)^m-\lambda_1\eta(t)>0$
and $\eta(t)>\eta_0$, for $t>0$. An integration of equation (\ref{3.15}) gives that
   \bess
T\leq\int_{\eta_0}^{\eta(T)}\frac{dr}{2r^{1+\frac{\alpha}{2}}+b^2q_1r^m-2\lambda_1r}
\leq\int_{\eta_0}^\infty\frac{dr}{2r^{1+\frac{\alpha}{2}}+b^2q_1r^m-2\lambda_1r}<\infty,
    \eess
which implies that $\eta(t)$ must blow up at a time $T^*\leq\int_{\eta_0}^\infty\frac{dr}{2r^{1+\frac{\alpha}{2}}+b^2q_1r^m-2\lambda_1r}$.
Hence this is a contradiction. This completes the proof. $\Box$

Before ending this section, we make the following remarks.

\begin{remark}\lbl{r3.2} $1.$ Theorem {\rm\ref{t3.2}} contains a new result.
First, we suppose there exists a positive constant $q_1$ such that $\inf_{x,y\in D}q(x,y)\geq q_1$.
When $b=0$, Fujita {\rm\cite{F}} showed that the solution of {\rm(\ref{3.9})} will exist
globally if the initial data is sufficiently small. Then we fixed the initial data sufficiently
small such that {\rm(\ref{3.9})} with $b=0$ has a unique global solution. Finally, we take the
suitable value of $b$ such that {\rm(\ref{3.10})} holds and it follows from
Theorem {\rm\ref{t3.2}} that the unique positive solution of {\rm(\ref{3.9})} will
blow up in finite time. Hence we can say that the noise induces the  finite time blow-up.

$2.$ From the proof of Theorem {\rm\ref{t3.2}}, we know that if one can prove that the solution of
{\rm(\ref{1.4})} is positive without using the property of $f(u)$, then the solution of {\rm(\ref{1.4})}
with $\sigma=u^m$ $(m>1)$ will
blow up in finite time under the condition that $f(u)\geq0$ for $u\geq0$.
Similar to   {\rm\cite{C2011}}, one can prove that Theorems
{\rm\ref{t3.1}} and {\rm\ref{t3.2}} also hold  for $ D=\mathbb{R}^n$ $($Theorem
{\rm 3.2} in {\rm\cite{C2011}}$)$.

$3.$ From {\rm(\ref{3.10})}, we   see  that for $m=1$ and $b^2q_1/2\geq\lambda_1$,
the solution of {\rm(\ref{3.9})} will blow up in finite time for any nonnegative initial data.
On the other hand, it follows from the proof of Theorem
{\rm\ref{t3.2}} that noise can make the existence time shorter.
 \end{remark}

\subsection{Third result on blow-up}
In this subsection, we consider the following equation
  \bes\left\{\begin{array}{lll}
   du=(\Delta u+f(u,x,t))dt+g(u,x,t)dW_t, \ \ \qquad & t>0,\ x\in  D,\\[1.5mm]
   u(x,0)=u_0(x), \ \ \ &\qquad\ \ x\in D,\\[1.5mm]
   u(x,t)=0, \ \ \ \ & t>0,\ x\in\partial D.
    \end{array}\right.\lbl{3.16}\ees
\begin{theo}\lbl{t3.3} Assume that all the assumptions in Theorem
{\rm\ref{t2.3}} hold. Assume further that $u_0$ is a nonnegative continuous function, $f(u,x,t)\geq0$
for $u\geq0,\,x\in D,\,t>0$ and
   \bes
&&g(u,x,t)\geq b^2u^m, \ \ \ m>1,\,b\in\mathbb{R},\nm\\
&&\inf_{x,y\in D}q(x,y)\geq q_1,\ \  \ \left(\int_ D u_0(x)\phi(x)dx\right)^{2(m-1)}\geq\frac{\lambda_1}{q_1b^2},
   \lbl{3.17}\ees
where $\lambda_1$ is defined as in {\rm(\ref{1.5})} and $q(x,y)$ is the
correlation function. Then the solution of {\rm(\ref{3.16})} will blow up in finite time
in $L^2(\Omega)\times L^\infty(D)$, that is,
there exists a constant $T^*>0$ such that
     \bes
\lim\limits_{t\rightarrow T^*-0}\left(\mathbb{E}\|u\|^2_{L^\infty(D)}\right)^{\frac{1}{2}}
=\lim\limits_{t\rightarrow T^*-0}\left(\mathbb{E}\sup_{x\in D}u(x,t)^2\right)^{\frac{1}{2}}=\infty.
    \lbl{3.18}\ees
\end{theo}

{\bf Proof.} By \cite{Cb,LR,PZ}, we know that (\ref{3.16}) has a unique local
solution. It follows from Theorem \ref{t2.3} that the solution of (\ref{3.16})
remains positive. Similar to the proof of Theorem \ref{t3.2}, it suffices to show that
$\mathbb{E}\hat u^2(t)$ blows up in finite time, where $\hat u(t)=(u,\phi)$.

By applying It\^{o}'s formula to $\hat u^2(t)$ and making use of (\ref{1.5}), we get
   \bes
\hat u^2(t)&=&(u_0,\phi)^2-2\lambda_1\int_0^t\hat u^2(s)ds+2
\int_0^t\int_ D\hat u(s) f(u,x,s)\phi(x)dxds\nm\\
&&+2\int_0^t\int_ D \hat u(s)g(u,x,s)\phi(x)dW(x,s)ds\nm\\
&&+\int_0^t\int_ D\int_ D q(x,y) g(u,x,s)\phi(x)g(u,y,s)\phi(y)dxdyds.
   \lbl{3.19}\ees
Let $\eta(t)=\mathbb{E}\hat u^2(t)$. By taking an expectation over (\ref{3.19}), we conclude that
  \bes
\eta(t)&=&(u_0,\phi)^2-2\lambda_1\int_0^t\eta(s)ds
+2\mathbb{E}\int_0^t\int_ D\hat u(s) f(u,x,s)\phi(x)dxds\nm\\
&&+\int_0^t\mathbb{E}\int_ D\int_ D q(x,y) g(u,x,s)\phi(x)g(u,y,s)\phi(y)dxdyds,
   \lbl{3.20}\ees
or, in the differential form
   \bes\left\{\begin{array}{lllll}
\ds\frac{d\eta(t)}{dt}=-2\lambda_1\eta(t)+2\mathbb{E}\hat u(t)\ds\int_ D f(u,x,t)\phi(x)dx\\[2.5mm]
\quad\quad\quad\ \ +\mathbb{E}\ds\int_ D\int_ D q(x,y) g(u,x,t)\phi(x)g(u,y,t)\phi(y)dxdy\\[2mm]
\eta(0)=\eta_0=(u_0,\phi)^2.
    \end{array}\right.\lbl{3.21}\ees
Again by Jensen's inequality, (\ref{3.21})
yields
     \bes\left\{\begin{array}{llll}
\ds\frac{d\eta(t)}{dt}\geq-2\lambda_1\eta(t)+q_1b^2\eta^m(t),  \\[1.5mm]
\eta(0)=\eta_0.
    \end{array}\right.\lbl{3.22}\ees
This implies that, for $\frac{b^2q_1}{2}\eta_0^m-\lambda_1\eta_0>0$,
we have $\frac{b^2q_1}{2}\eta(t)^m-\lambda_1\eta(t)>0$
and $\eta(t)>\eta_0$ for $t>0$. An integration of equation (\ref{3.22}) gives that
   \bess
T\leq\int_{\eta_0}^{\eta(T)}\frac{dr}{b^2q_1r^m-2\lambda_1r}
\leq\int_{\eta_0}^\infty\frac{dr}{b^2q_1r^m-2\lambda_1r}<\infty,
    \eess
which implies $\eta(t)$ must blow up at a time $T^*\leq\int_{\mu_0}^\infty\frac{dr}{b^2q_1r^m-2\lambda_1r}$.
Hence this is a contradiction. This completes the proof. $\Box$

\begin{remark}\lbl{r3.3}   Theorem {\rm\ref{t3.3}} holds for
$g(u,x,t)=bu^m$, $m=2,3,\cdots$. It shows that the noise can induce a
singularity.
\end{remark}

\section{Global solution for a stochastic Allen-Cahn equation}
\setcounter{equation}{0}

In this section, we show that the solution of a  stochastic Allen-Cahn equation does not have finite time singularities and it exists globally. This is an example of SPDEs whose coefficients are not globally Lipschitz continuous.

%consider the global solution to stochastic parabolic It\^{o} equations. The earlier results about the
%global solution to stochastic partial differential equations often
%require the stochastic term satisfying the global Lipschitz condition \cite{DL,PZ1992}. In this section, we give another example: the

We consider the following stochastic Allen-Cahn equation.
      \bes\left\{\begin{array}{lll}
   du=(\Delta u+u(1-u^2))dt+bu^mdW_t, \ \ \qquad &t>0,\ x\in  D,\\[1.5mm]
   u(x,0)=u_0(x), \ \ \ &\qquad\ \ x\in D,\\[1.5mm]
   u(x,t)=0, \ \ \ \ &t>0,\ x\in\partial D,
    \end{array}\right.\lbl{4.1}\ees
where $1<m<2,\,b\in\mathbb{R}$. If $b=0$, (\ref{4.1}) becomes the well-known
deterministic Allen-Cahn equation \cite{AC},
which describes the process of phase separation in iron alloys,
including order-disorder transitions. Hairer et al. \cite{HRW} considered (\ref{4.1}) with
$m=0$.   The equation (\ref{4.1}) with $b=0$ has a global solution. We
want to know when the solution of (\ref{4.1}) exists globally and when the solution
blows up. In this section, we shall
use the Lyapunov functional method to prove that the solution of (\ref{4.1}) exists globally, i.e., no finite time blow up.

Throughout this section, we assume that $H=L^2( D)$, $\|\cdot\|=\|\cdot\|_{H}$. Moreover,
$W(x,t)$ is a Wiener random field, and $q(x,y)$ is its covariance function as defined   in
Section 2.

Let $u(x,t;u_0)$ be a solution of (\ref{4.1}) with the initial data $u_0$.
We first give the definition of global solution.
\begin{defi} \lbl{4.1} A function $u(x,t)\in H( D)\cap H_0^1( D)$ is said to be non-explosive solution of
{\rm(\ref{4.1})} if
    \bess
\lim\limits_{r\rightarrow\infty}P\{\sup_{0\leq t\leq T}\|u_t\|>r\}=0,
   \eess
for any $T>0$. If the above holds for $T=\infty$, the solution $u(x,t)$ is said to
be ultimately bounded, i.e., global solution.
\end{defi}

We shall use Lyapunov functional method to obtain the existence of global solution
to (\ref{4.1}). In the following, we recall the definition of Lyapunov functional (\cite{Cb}).
We do this for a  more general stochastic partial differential equation
   \bes\left\{\begin{array}{lll}
du=(A u+F(u))dt+{\sum}(u)dW_t, \ \ \ t\geq0,\\
u(x,0)=h(x),
   \end{array}\right.
   \lbl{4.2}\ees
where $A,\,F$ and $\sum$ are assumed to be non-random or deterministic.
Let $V$ be a separate Hilbert space.
Here we say that a $\mathcal {F}_t$-adapted $V$-valued process $u$ is a
\emph{strong solution} of equation (\ref{4.2}) if $u\in L^2(\Omega\times[0,T];V)$,
and for any $\phi\in V$, the following equation
   \bess
(u,\phi)=(h,\phi)+\int_0^t\langle Au+F(u),\phi\rangle ds+\int_0^t(\phi,\sum(u)dW_s)
   \eess
holds for each $t\in[0,T]$ a.s.
Recall that the generator for this stochastic partial differential equation is (see Chapter 7 in \cite{Cb})
   \bes
\mathcal {L}_t\Phi(v,t)&=&\frac{\partial}{\partial_s}\Phi(v,s)+\frac{1}{2}Tr[\Phi''(v,t){\sum}_t(v)Q{\sum}^*_t(v)]\nm\\
&&+\langle A_tv,\Phi'(v,t)\rangle+(F_t(v),\Phi'(v,t)),
   \lbl{4.3}\ees
where $Q$ is covariance operator. Let $U\subset H$ be a neighborhood of the origin. A function
$\Phi$: $U\times\mathbb{R}^+\rightarrow\mathbb{R}$ is said to be a \emph{Lyapunov functional}
for the equation (\ref{4.2}), if

(1) $\Phi$ is locally bounded and continuous such that its first two partial derivatives
$\partial_t\Phi(v,t),\,\partial_x\Phi(v,t)$ and $\partial_{xx}\Phi(v,t)$ exist, and
$\partial_t\Phi(v,t),\,\partial_x\Phi(v,t)$ are locally bounded.

(2) $\Phi(0,t)=0$ for all $t\geq0$, and, for any $r>0$, there is $\delta>0$ such that
  \bess
\inf_{t\geq0,\|v\|\geq r}\Phi(v,t)\geq\delta.
   \eess

(3) For every  $t\geq0$ and $v\in U\cap H^1_0$,
   \bess
\mathcal {L}_t\Phi(v,t)\leq0,
  \eess
where $\partial_x$ and $\partial_{xx}$ are Fr\'{e}chet derivative, see \cite[p.196-201]{Cb}.

Let $U\times\mathbb{R}^+\rightarrow\mathbb{R}^+$ be a Lyapunov
functional and let $u_t$ denote the strong solution of {\rm\ref{4.2}} with initial data
$u_0$. For $r>0$, let $B_r=\{h\in H:\|h\|\leq r\}$ such that $B_r\subset U$. Define
   \bess
\tau=\inf\{t>0:\ u_t\in B_r^c,\ u_0\in B_r\},
   \eess
with $B_r^c=H\ B_r$. We put $\tau=T$ if the set is empty. Then the
process $\phi_t=\Phi(u_{t\wedge\tau},t\wedge\tau)$ is a local $\mathcal {F}_t$-supermartingale
and the following Chebyshev inequality holds
   \bess
P\{\sup_{0\leq t\leq T}\|u_t\|>r\}\leq\frac{\Phi(u_0,0)}{\Phi_r},
  \eess
where
   \bess
\Phi_r=\inf_{0\leq t\leq T,h\in U\cap B_r^c}\Phi(h,t).
  \eess

In order to obtain the global solution of (\ref{4.1}), we need
the following Lemma.

  \begin{lem}\lbl{l4.1} {\rm\cite[p.200, Theorem 3.2]{Cb}} Suppose that
there exists a Lyapunov functional $\Phi:\ H\times\mathbb{R}_+\rightarrow\mathbb{R}_+$ such
that
   \bess
\Phi_r=\inf_{t\geq0,\|h\|\geq r}\Phi(h,t)\rightarrow\infty,\ \ {\rm as}\ r\rightarrow\infty.
   \eess
Then the solution $u_t$ is ultimately bounded.
\end{lem}

Let $\Phi(u,t)=e^{-\alpha t}\Psi(u,t)$. Then it follows from Lemma \ref{l4.1} that
if there exists a
constant $\alpha>0$ such that
   \bess
\mathcal {L}_t\Psi\leq\alpha\Psi(v,t)\ \ \ \ {\rm for}\ {\rm any}\ v\in H_0^1,
   \eess
and the infimum $\inf_{t\geq0,\|v\|\geq r}\Psi(v,t)=\Psi_r$ exists such that
$\lim\limits_{r\rightarrow\infty}\Psi_r=\infty$, then the solution $u$ does not
explode in finite time.

Now, we use the Lemma \ref{l4.1} to examine the global solution of (\ref{4.1}).

\begin{theo}\lbl{t4.1} Assume that $1<m<2$ and $u_0(x)\geq0$ for $x\in\bar D$.
Assume further that there exists a positive constant $q_0$ such that the covariance function $q(x,y)$ satisfies the condition $\sup_{x,y\in\bar D}q(x,y)\leq q_0$.
Then {\rm(\ref{4.1})} has a strong global solution.
   \end{theo}

{\bf Proof.} It follows from \cite{Cb,LR,PZ} that (\ref{4.1}) has a local solution
on $[0,T]\times D$. By Theorem \ref{t2.1},   this local solution
is positive. Now, we use Lemma \ref{l4.1} to prove the solution does not blow
up in finite time. Define $\Psi(v,t)=\|v\|^2$. Then $\lim\limits_{r\rightarrow\infty}\inf_{t\geq0,\|v\|\geq r}\Psi(v,t)=\infty$.
Direct calculations show that
  \bes
\mathcal {L}_t\Psi(v,t)&=&\frac{\partial}{\partial_s}\Psi(v,s)+\frac{1}{2}Tr[\Psi''(v,t)v_t^mQv^m_t]\nm\\
&&+\langle \Delta v,\Psi'(v,t)\rangle+(v-v^3,\Psi'(v,t))\nm\\
&=&-2\int_ D |\nabla v|^2dx+2\int_ D(v^2-v^4)dx+\int_ D q(x,x)v^{2m}(x)dx\nm\\
&\leq&2\|v\|^2-2\|v\|_{L^4}^4+q_0\|v\|_{L^{2m}}^{2m}.
   \lbl{4.4}\ees
By using interpolation inequality (\ref{2.4}) with $r=2m,\,p=2$ and $q=4$, we have
  \bes
q_0\|v\|_{L^{2m}}^{2m}&\leq&C\|v\|_{L^2}^{2m\theta}\|v\|_{L^{4}}^{2m(1-\theta)}\nm\\
&\leq&\epsilon\|v\|_{L^{4}}^{2m(1-\theta)\frac{1}{1-m\theta}}+C(\epsilon)\|v\|_{L^2}^{2}\nm\\
&=&\epsilon\|v\|_{L^{4}}^4+C(\epsilon)\|v\|_{L^2}^{2},
   \lbl{4.5}\ees
where $\theta=\frac{2-m}{m}$ satisfying (\ref{2.5}). Substituting (\ref{4.5}) into (\ref{4.4}), we have
    \bess
\mathcal {L}_t\Psi(v,t)\leq C\|v\|^2-\|v\|_{L^4}^4\leq C\|v\|^2=C\Psi(v,t),
    \eess
which implies that all the assumptions in Lemma \ref{l4.1} hold. Thus by Lemma \ref{l4.1}
we know that the solution of (\ref{4.1}) exists globally. This completes the proof. $\Box$

It can be shown that Theorem \ref{t4.1} also holds if (\ref{4.1}) is replaced by
      \bes\left\{\begin{array}{lll}
   du=(\Delta u-u^\gamma)dt+bu^mdW_t, \ \ \qquad &t>0,\ x\in  D,\\[1.5mm]
   u(x,0)=u_0(x), \ \ \ &\qquad\ \ x\in D,\\[1.5mm]
   u(x,t)=0, \ \ \ \ &t>0,\ x\in\partial D,
    \end{array}\right.\lbl{4.6}\ees
for  $b\in\mathbb{R},\,1<m<(\gamma+1)/2$ and $\gamma>1$ satisfying $(-1)^\gamma=-1$.

\begin{col}\lbl{c4.1}
Assume that $\gamma>1$, $1<m<(\gamma+1)/2$,
$u_0\geq0$ and assume also that there exists a positive constant $q_0$ such that the covariance function $q(x,y)$ satisfies the condition $\sup_{x,y\in\bar D}q(x,y)\leq q_0$.
Then {\rm(\ref{4.6})} has a unique strong global solution.
\end{col}

Theorem \ref{t4.1} and Corollary \ref{c4.1} imply that if the nonlinearity $f(u)=ku-u^\gamma$ can
control the stochastic term $u^m$, i.e., $m<(\gamma+1)/2$, the stochastic partial differential
equation also has global solutions, which is different from the earlier results.

Similarly, we can use Lemma \ref{l4.1} to study the following
stochastic partial differential equation
      \bes\left\{\begin{array}{lll}
   du=(\nu\Delta u+au(1-u^2))dt+k\ds\sum_{i=1}^3\ds\frac{\partial u}{\partial_{x_i}}dW_i(x,t), \ \ \qquad &t>0,\ x\in  D,\\[2.5mm]
   u(x,0)=u_0(x), \ \ \ &\qquad\ \ x\in D,\\[1.5mm]
   u(x,t)=0, \ \ \ \ &t>0,\ x\in\partial D,
    \end{array}\right.\lbl{4.7}\ees
where $ D\subset\mathbb{R}^3$ is a bounded domain with smooth boundary
$\partial D$ and $W_i(x,t)$ are Wiener random fields with bounded covariance functions $q_{jk}(x,y)$ such that
   \bess
\sum_{j,k=1}^3q_{jk}(x,x)\xi_j\xi_k\leq q_0|\xi|^2,\ \ \ \forall\xi\in\mathbb{R}^3
   \eess
for some $q_0>0$. From \cite{Cb,LR,PZ}, we know that equation (\ref{4.7}) has a strong
solution $u\in H_0^1$ (see Theorem 6-7.5 in \cite{Cb}). Similar to the proof of Theorem \ref{t4.1}, we
define $\Psi(v,t)=\|v\|^2$. Then $\lim\limits_{r\rightarrow\infty}\inf_{t\geq0,\|v\|\geq r}\Psi(v,t)=\infty$.
Again, we have
  \bess
\mathcal {L}_t\Psi(v,t)&=&-2\nu\int_ D |\nabla v|^2dx+2a\int_ D(v^2-v^4)dx\nm\\[2mm]
&&+\int_ D \sum_{j,k=1}^3q_{jk}(x,x)\frac{\partial v(x)}{\partial_{x_j}}\frac{\partial v(x)}{\partial_{x_k}}dx\nm\\[2mm]
&\leq&2a\|v\|^2-(2\nu-q_0)\int_ D |\nabla v|^2dx\\
&\leq&2a\|v\|^2
   \eess
provided that $2\nu-q_0>0$. Using Lemmas \ref{l4.1}, we
obtain the following result.
\begin{theo}\lbl{t4.2} Assume that $2\nu-q_0>0$ and $u_0(x)\geq0$ for $x\in\bar D$.
Then {\rm(\ref{4.7})} has a unique strong global solution.
   \end{theo}

%Similarly, we can deal with the case that $f(u)=u-u^\gamma$ and $\sigma(u,\nabla u,x,t)=\nabla u$ in (\ref{1.4}).

\medskip

\noindent {\bf Acknowledgment} The first author was supported in part
by NSFC of China grants 11301146, 11171064 and 11226168. The second author was partially supported
by the NSF grant 1025422. Part of this work was done while Guangying Lv was visiting Illinois  Institute of Technology, Chicago,   USA.
The authors are grateful
to the referees for their valuable suggestions and comments on the
original manuscript.

 \end{document}